\RequirePackage{amsmath}
\documentclass[runningheads]{llncs}
\usepackage{graphicx}
\usepackage{amssymb}
\usepackage{mathtools,xfrac}
\usepackage{url}
\usepackage{color}
\usepackage{tikz}
\usepackage[utf8]{inputenc}
\usepackage{hyperref}

\newcommand{\Z}{\mathbb{Z}}
\newcommand{\R}{\mathbb{R}}
\newcommand{\qedEQN}{\tag*{$\qed$}}

\DeclareMathOperator{\supp}{supp}

\DeclareMathOperator{\rank}{rank}

\newtheorem{thm}{Theorem}
\newtheorem{lem}[thm]{Lemma}
\newtheorem{cor}[thm]{Corollary}

\newtheorem{CLAIM}[thm]{Claim}
\makeatletter
\newcommand*{\transpose}{%
  {\mathpalette\@transpose{}}%
}
\newcommand*{\@transpose}[2]{%
  \raisebox{\depth}{\m@th$#1\intercal$}%
}
\makeatother

\begin{document}

\title{Improving proximity bounds using sparsity}
%
%
\author{Jon Lee\inst{1} \and
Joseph Paat\inst{2} \and
Ingo Stallknecht\inst{2} \and
Luze Xu\inst{1}}
\authorrunning{Lee, Paat, Stallknecht \& Xu}
%
\institute{Department of Industrial and Operations Engineering, University of Michigan, USA \and
Department of Mathematics, ETH Z\"{u}rich, Switzerland}
\maketitle              
\begin{abstract}
We refer to the distance between optimal solutions of integer programs and their linear relaxations as \emph{proximity}.
In 2018, Eisenbrand and Weismantel proved that proximity is independent of the dimension for programs in standard form.
We improve their bounds using existing and novel results on the \emph{sparsity} of integer solutions.
We first bound proximity in terms of the largest absolute value of any full-dimensional minor in the constraint matrix, and this bound is tight up to a polynomial factor in the number of constraints.
We also give an improved bound in terms of the largest absolute entry in the constraint matrix, after efficiently transforming the program into an equivalent one.
Our results are stated in terms of general sparsity bounds, so any new results on sparse solutions immediately improves our work.
Generalizations to mixed integer programs are also discussed.
\end{abstract}
%
%
\section{Introduction}

Let $A \in \Z^{m\times n}$ with $\rank(A) = m$, $c \in \Z^n$, and $b \in \Z^m$.
For $k \in \{1, \dotsc, m\}$, denote the largest absolute value of a minor of order $k$ of the matrix $A$ by
\[
\Delta_k := \Delta_k(A) := \max \{|\det(B)| : B \text{ is a } k \times k \text{ submatrix of } A\}.
\]
Note that $\Delta_1 = \|A\|_{\infty}$ is the largest absolute entry of $A$. For simplicity, we set $\Delta:=\Delta_m$.

We consider the standard form integer program
\begin{equation}\label{IP}
\max\{c^\transpose z : Az = b, ~  z \in \Z^n_{\ge 0} \} \tag{\text{IP}}
\end{equation}
and its linear relaxation
\begin{equation}\label{LP}
\max \{ c^\transpose x :  Ax =  b, ~ x \in \R^n_{\ge 0}\} \tag{\text{LP}}.
\end{equation}
Throughout, we assume that~\eqref{IP} is both feasible and bounded.

Given an optimal vertex solution $x^*$ to~\eqref{LP}, we investigate the question of~\eqref{LP} to~\eqref{IP} proximity: can we bound the distance from $x^*$ to some optimal solution $z^*$ of~\eqref{IP}?
We refer to any bound $\tau$ on $\|x^* - z^*\|_1$ as a \emph{proximity bound}.
Proximity bounds have a variety of implications in the theory of integer programming.
For example, a proximity bound of $\tau$ translates into a bound of $\tau \cdot \|c\|_{\infty}$ on the so-called integrality gap~\cite{CooSens,GraProx,HochConv}.
Furthermore, strong proximity bounds reduce the time needed for a local search algorithm to find an optimal~\eqref{IP} solution starting from an optimal~\eqref{LP} solution~\cite{WeiProx}.
%
%

One of the first seminal results on proximity is by Cook et al.~\cite{CooSens}, who established that there exists an optimal solution $z^*$ to~\eqref{IP} satisfying

\begin{equation}\label{eqCookBound_infty}
\|z^* - x^*\|_\infty \le n \cdot \max\{\Delta_{k} : k \in \{1, \dotsc, m\}\} .
\end{equation}
Cook et al. actually consider problems in inequality form, i.e., with constraints $Ax \le b$ rather than $Ax = b$, but their results easily translate to the standard form setting.
 A closer analysis reveals that $\Delta$ suffices for the standard form problem rather than $\max\{\Delta_{k} : k \in \{1, \dotsc, m\}\} $ stated in~\eqref{eqCookBound_infty}.
Furthermore, if we naively extend~\eqref{eqCookBound_infty} to a bound on proximity in terms of the $\ell_1$-norm, then we obtain $\|z^* - x^*\|_1 \le n^2\Delta$.
Another closer analysis
gives us the bound
\begin{equation}\label{eqCookBound}
\|z^* - x^*\|_1 \le (m+1)n\Delta.
\end{equation}
See the proof of Lemma~\ref{MixedProx} for the two `closer analyses' referred to above.
Cook et al.'s bound has been generalized to various problems including those with separable convex objective functions~\cite{GraProx,HochConv,WerSep} or with mixed integer constraints~\cite{paat2018distances}, and extended to alternative data parameters such as $k$-regularity~\cite{Lee_LAAsubspaces,LeeRegu} and the magnitude of Graver basis elements~\cite{eisenbrand2019algorithmic}.

Notice that the proximity bound in~\eqref{eqCookBound} depends on the dimension.
In 2018, Eisenbrand and Weismantel~\cite{WeiProx} proved that proximity can be bounded independently of the dimension by establishing the bound
\begin{equation}\label{eqEWEntryBound}
\|z^* - x^*\|_1 \le m(2m\cdot \|A\|_{\infty}+1)^m.
\end{equation}
The proof of~\eqref{eqEWEntryBound} involves the Steinitz Lemma~\cite{S1913}.
By applying a different norm in the Steinitz Lemma, Eisenbrand and Weismantel's result can be restated as
\begin{equation}\label{eqEWDeltaBound}
\|z^* - x^*\|_1 \le m(2m+1)^m\cdot \Delta.
\end{equation}
The proximity bounds~\eqref{eqEWEntryBound} and~\eqref{eqEWDeltaBound} also hold for standard form problems with additional upper bound constraints on the variables.

As for lower bounds on proximity, it is not difficult to come up with examples demonstrating $\|z^* - x^*\|_1 \ge m \cdot \Delta$ and $\|z^* - x^*\|_1 \ge \|A\|_{\infty}^m$.
Aliev et al.~\cite{AHO2019} give a tight lower bound $\Delta-1$ on proximity in terms of the $\ell_\infty$-norm when $m = 1$.
However, it remains an open question if \eqref{eqCookBound}, \eqref{eqEWEntryBound}, or \eqref{eqEWDeltaBound} is tight in general.
To this end, Oertel et al. established that the upper bound $\|z^* - x^*\|_1 \le (m+1) \cdot (\Delta-1)$ holds for most problems, where `most' is defined parametrically with $b$ treated as input~\cite{OPW2020}.

\subsection{Statement of Results and Overview of Proof Techniques}

The focus of this paper is to create stronger proximity bounds.
Recall that $\Delta := \Delta_m$.
Our first main result is an improvement over~\eqref{eqEWDeltaBound} for fixed $\Delta$.
We always consider the logarithm $\log(\cdot)$ to have base two.
\begin{thm}\label{DeltaProx}
For every optimal \eqref{LP} vertex solution $x^*$, there exists an optimal \eqref{IP} solution $z^*$ such that
\[
\|z^* -x^*\|_1 <3m^2 \log (2\sqrt{m}\cdot \Delta^{1/m}) \cdot \Delta.
\]
\end{thm}
Theorem~\ref{DeltaProx} demonstrates that proximity in the $\ell_1$-norm between~\eqref{IP} solutions and~\eqref{LP} \emph{vertex} solutions is bounded by a polynomial in $m$ and $\Delta$.
%
We focus on vertex solutions because proximity may depend on $n$ for general non-vertex~\eqref{LP} solutions.
For example, suppose that $c = 0$ and take any feasible solution to~\eqref{LP}, which is optimal in this case, such that each of the $n$ components is in $\sfrac12+\Z$.

Our second main proximity result is in terms of $\|A\|_{\infty}$ \emph{after $A$ is transformed by a suitable unimodular matrix}.
Recall that a unimodular matrix $U \in \Z^{m\times m}$ satisfies $|\det(U)| = 1$, so the $m\times m$ minors of $UA$ have the same magnitudes as those of $A$.
Moreover, the optimal solutions of~\eqref{IP} are the same as the optimal solutions to
\begin{equation}\label{eqUIP}
\max\{c^\transpose z : UAz = Ub, ~ z \in \Z^n_{\ge 0}\}.\tag{$U$ - IP}
\end{equation}
%
%
If we know an $m \times m$ submatrix $B$ of $A$ satisfying $|\det(B)| = \Delta$, then in polynomial time we can find a unimodular matrix $U$ such that $UB$ is upper triangular, which implies that $\Delta \leq \|UB\|_\infty ^m$.
Here, we can apply Theorem~\ref{DeltaProx} to obtain the bound
\begin{equation}\label{eqLargestHNF}
\|z^* -x^*\|_1 <3m^2 \log (2\sqrt{m}  \cdot \|UB\|_\infty )\cdot \|UB\|_\infty ^m.
\end{equation}

The previous bound is predicated on the knowledge of an $m \times m$ submatrix of maximum absolute determinant, which is NP-hard to find~\cite{K1995}.
However, Di Summa et al.~\cite{SumSubdet} established that this submatrix can be approximated in polynomial time.
In particular, they demonstrated that there exists an $m \times m$ submatrix $B$ of $A$ satisfying
\begin{equation}\label{eqDetApprox}
\Delta \le |\det(B)| \cdot (2\log(m+1))^{m/2}
\end{equation}
that can be found in time polynomial in $m$ and $n$.\footnote{The approximation result of Di Summa et al. involves an $\epsilon$ factor of precision and the running time is polynomial in $m,n,1/\epsilon$. For the sake of presentation, we have fixed this $\epsilon$ to $1/m$ and obtain a polynomial time algorithm in $m,n$.}
Similarly to the above, there exists a unimodular matrix $U \in \Z^{m\times m}$ for which $UB$ is upper triangular that can also be found in polynomial time.
We can use this approximate largest absolute determinant to derive our second main result.
We denote the linear relaxation of~\eqref{eqUIP} by ($U$ - LP).
\begin{thm}\label{eqThmEntry}
Let $B$ be an $m\times m$ submatrix $B$ of $A$ satisfying~\eqref{eqDetApprox} and $U \in \Z^{m\times m}$ a unimodular matrix such that $UB$ is upper triangular.
Then for every optimal \textnormal{($U$ - LP)} vertex solution $x^*$, there exists an optimal \eqref{eqUIP} solution $z^*$ satisfying
\[
\|z^* -x^*\|_1 <3m^2 \log (2\sqrt{2m\log (m+1)} \cdot \|UB\|_\infty  )  \cdot (2\log (m+1))^{m/2}\cdot \|UB\|_\infty ^m.
\]
\end{thm}
%
%
It is worth reemphasizing that the proximity bound in Theorem~\ref{eqThmEntry} can be determined in polynomial time, which is in contrast to the bound in~\eqref{eqLargestHNF}, and the dependence on $m$ is significantly less than the bound in~\eqref{eqEWEntryBound}.

The proofs of Theorems~\ref{DeltaProx} and~\ref{eqThmEntry} are based on combining proof techniques of Cook et al.~\cite{CooSens} with results on the sparsity of optimal solutions to~\eqref{IP}.
The \emph{support} of a vector $x \in \R^n$ is defined as
\[
\supp(x) := \{i \in \{1, \dotsc, n\} : x_i \neq 0\}.
\]
A classic theorem of Carath\'eodory states that $|\supp(x^*)| \le m$ for every optimal vertex solution of~\eqref{LP}.
It turns out that the minimum support of an optimal solution to~\eqref{IP} is not much larger.
Denote the minimum support of an optimal solution of~\eqref{IP} by
%
\[
S := \min\left\{|\supp(z^*)| :
z^* \text{ is an optimal solution for~\eqref{IP}}
\right\}.
\]
Aliev et al.~\cite{ADEOW2018} established that
\begin{equation}\label{eqSparseSolution}
S \le m+\log\left(\sqrt{\det(AA\strut^{\transpose})}\right) \le  2m\log(2\sqrt{m} \cdot \|A\|_{\infty}).
\end{equation}
For other results regarding sparsity, see~\cite{ADOO2017,ES2006,OPW2019} for general $A$ and~\cite{BG2004,BrunsGHMW1999,CookFS1986,Sebo1990} for matrices that form a Hilbert basis.
See also the manuscript of Aliev et al.~\cite{AADO2019}, who give improved sparsity bounds for \emph{feasible solutions} to special classes of integer programs and provide efficient algorithms for finding such solutions.
Using sparsity, we derive the following proximity bound which forms the basis for Theorems~\ref{DeltaProx} and~\ref{eqThmEntry}.

\begin{lem}\label{MixedProx}
For every optimal \eqref{LP} vertex solution $x^*$, there exists an optimal \eqref{IP} solution $z^*$ such that
\[
\| z^* -x^* \|_1 < (m+1)\cdot(m+S)\cdot \Delta.
\]
\end{lem}

Lemma~\ref{MixedProx} improves~\eqref{eqCookBound} by replacing the dependence on the dimension $n$ to $m+S$.
Lemma~\ref{MixedProx} is stated for a generic sparsity bound, so one could use it together with~\eqref{eqSparseSolution} to achieve a proximity bound in terms of $\Delta$ and $\|A\|_{\infty}$.
In order to provide a bound for proximity that is uniform in the data parameter, we prove a new sparsity result in terms of $\Delta$.
\begin{thm}\label{DeltaSpars}
There exists an optimal \eqref{IP} solution $z^*$ such that
\begin{equation*}
|\supp (z^*)|< 2m\log \bigl(\sqrt{2m} \cdot \Delta^{1/m}\bigr).
\end{equation*}

\end{thm}

%
%
\smallskip

Our proximity bounds can be generalized to mixed integer programs.
Given an index set $\mathcal{I} \subseteq \{1, \dotsc, n\}$, the mixed integer program with integrality constraints indexed by $\mathcal{I}$ is
\begin{equation}\label{MIP}
    \max\{c^\transpose z : Az = b, ~  z\ge0, ~ z_i \in \Z ~\forall~i\in\mathcal{I} \}. \tag{\text{MIP}}
\end{equation}
%
Similarly to~\cite[Corollary 4]{ADEOW2018}, we establish the extension of Theorem \ref{DeltaSpars} to \eqref{MIP}.
\begin{cor}\label{DeltaSparsMIP}
    There exists an optimal \eqref{MIP} solution $z^*$ satisfying
    \begin{equation*}
    |\supp (z^*)|< m+2m\log \bigl(\sqrt{2m}\cdot \Delta^{1/m}\bigr)=2m\log \bigl(2\sqrt{m}\cdot \Delta^{1/m}\bigr).
    \end{equation*}

\end{cor}

We obtain the following proximity result by applying Corollary \ref{DeltaSparsMIP}.
\begin{cor}\label{DeltaProxMIP}
    For every optimal \eqref{LP} vertex solution $x^*$, there exists an optimal \eqref{MIP} solution $z^*$ such that
    \[
    \|z^* -x^*\|_1 <3m^2 \log (2\sqrt{2m}\cdot \Delta^{1/m}) \cdot \Delta.
    \]
\end{cor}

Our results also extend to integer programs in general form.
Let $A\in\mathbb{Z}^{m \times n}$ and $B\in\mathbb{Z}^{m \times d}$ be matrices satisfying $\rank ([A, B])=m$.
Note that it is not necessary to assume that $\rank(A) = m$ in our general form results.
Let $C\in\mathbb{Z}^{t \times d}$, $c \in \Z^{n+d}$, $b_1 \in \Z^m$, and $b_2 \in \Z^{t}$.
The \emph{general form integer program} is
\begin{equation}\label{GIP}
\max \left\{c^\transpose z  :
\begin{array}{rcl}
[A,~ B]~z  & = &  b_1\\[.15 cm]
[{0},~C]~z  & \le & b_2
\end{array},
~~z \in \Z^{n+d},
~~z_i \ge 0~\forall ~i \in \{1, \dotsc, n\}
\right\}.
\tag{GIP}
\end{equation}
We define the general form linear program~(GLP) similarly.
Previously cited bounds on proximity hold for~\eqref{GIP}.
%
However, our analysis reveals that proximity for~\eqref{GIP} depends on the potentially smaller data parameter
\[
\delta := \max \left\{|\det(E)| :
\begin{array}{l}
E \text{ is any submatrix of } \begin{pmatrix} A & B \\ 0 & C\end{pmatrix} \\[.25 cm]
\text{ defined using the first } m \text{ rows}
\end{array}
\right\}.
\]
If $t =0$ and $d=0$, then~\eqref{GIP} is a standard form problem and $\delta = \Delta_m(A)$.
If $m =0$ and $n=0$, then~\eqref{GIP} is an inequality form problem and $\delta = \max\{\Delta_k(C) : k \in \{1, \dotsc, d\}\}$.
\begin{cor}\label{corProxPlusSparsity}
For every optimal \textnormal{(GLP)} vertex solution $x^*$, there exists an optimal \eqref{GIP} solution $z^*$ such that
\[
\|z^* -x^* \|_1 < \min \{ m+t+1,n+d \}\cdot  \Bigl(\min \bigl\{ n,2m \cdot \log (2\sqrt{m} \cdot \delta^{1/m}) \bigr\} +d\Bigr)\cdot  \delta .
\]
\end{cor}
The proximity bound in Corollary~\ref{corProxPlusSparsity} coincides with the best known bounds in both the standard setting and the inequality form setting.

Going beyond integer linear optimization problems, it would be ideal for proximity bounds in terms of sparsity to extend to integer programs with separable convex objective functions (see~\cite{GraProx,HochConv} for similarities between the linear and separable convex setting).
However, for separable convex \emph{maximization problems}, strong proximity bounds do not exist for exact solutions, in general, even though sparsity results apply.
In contrast, for separable convex \emph{minimization problems}, strong sparsity bounds do not exist for exact solutions, in general, even though the classic proximity techniques apply.


In the remainder of this paper, we provide the proofs of our results. In Section \ref{sec:sparsity}, we first present our proofs of the proximity bound derived from a generic sparsity bound (Lemma~\ref{MixedProx}) and the novel sparsity bound (Theorem~\ref{DeltaSpars} and Corollary~\ref{DeltaSparsMIP}), because the proofs of the proximity results depend on these.
Then in Section \ref{sec:proximity},  we
provide the proofs of the proximity results for the standard form integer programs (Theorem~\ref{DeltaProx} and Theorem~\ref{eqThmEntry}). The proof of the mixed integer case (Corollary~\ref{DeltaProxMIP}) is omitted because it is the same as the proof of the pure integer case except that a different sparsity bound is applied. Additionally, we provide a proof of the proximity result in the general form setting (Corollary~\ref{corProxPlusSparsity}).

%
\section{Proofs regarding sparsity}\label{sec:sparsity}

Given $A \in \R^{m \times n}$ and $I \subseteq \{1, \dotsc, n\}$, we let $A_I \in \R^{m \times |I|}$ denote the columns of $A$ indexed by $I$.
If $I = \{i\}$ for some $i \in \{1, \dotsc, n\}$, then $A_i := A_{I}$.
Similarly, given $u \in \R^n$, we let $u_I \in \R^{|I|}$ denote the components of $u$ indexed by $I$.

\begin{proof}[of Lemma~\ref{MixedProx}]
We prove the result by projecting the optimization problems onto the union of the supports of $x^*$ and an optimal~\eqref{IP} solution with minimal support.
 Let $\bar{z}  \in\mathbb{Z}^{n}_{\ge 0}$ be an optimal \eqref{IP} solution with minimum support.
 By the definition of $S$ we have $|\supp (\bar{z})|= S$.
 As $x^*$ is an optimal vertex solution of \eqref{LP}, we also have $|\supp (x^*)| \le m$.
 Define
 \[
 H :=  \supp (x^*) \cup \supp (\bar{z}),
 \]
 and note that
 \begin{equation}\label{eqDefnH}
 |H| = |\supp(x^*) \cup \supp (\bar{z})|\le |\supp(x^*) | +  |\supp(\bar{z}) | \leq m+S.
 \end{equation}
If $n=m$, then $A$ is invertible and there exists a unique solution $A^{-1}b$ to the system $Ax = b$.
In this case $x^* = \bar{z} = A^{-1}b$.
%
Therefore, $\|x^* - \bar{z}\|_1 = 0$.
For the rest of the proof, we assume that $n>m$ and $H=\{ 1,\dots ,|H| \}$.

Consider the optimization problems

\begin{equation}\label{IP2}
\max \big\{ c_H^\transpose z : A_H z   =  b ,~ z   \in \Z^{|H|}_{\ge 0}\big\} \tag{IP2}
\end{equation}
and
\begin{equation}\label{LP2}
\max \big\{ c_H^\transpose x : A_Hx    =  b,~ x   \in  \R^{|H|}_{\ge 0} \big\}. \tag{LP2}
\end{equation}
%
%
Observe that $x^*_H$ is an optimal basic feasible solution for~\eqref{LP2}, and $\bar{z}_H$ is an optimal solution for~\eqref{IP2}.

Rewrite \eqref{LP2} in inequality form as follows:
\[
\max\big\{ c_H^\transpose x   : A_Hx =b, -I_Hx \le 0,~x \in\mathbb{R}^{|H| } \big\},
\]
where $I_{H}$ is the $|H|\times |H|$ identity matrix.
Partition the rows of $-I_H$ into $D_1$ and $D_2$ such that $D_1 \bar{z}_H < D_1 x^*_H  $ and $D_2 \bar{z}_H \geq D_2 x^*_H $.
Define the polyhedral cone
\begin{equation}\label{eqCone}
K:= \left\{ u\in\mathbb{R}^{|H| }  :  A_H u = 0, ~ D_1 u \le 0, ~D_2 u \ge 0\right\}.
\end{equation}
Observe that $\bar{z}_H - x^*_H \in K$.

The cone $K$ is pointed because the submatrix defined by $D_1$ and $D_2$ in~\eqref{eqCone} has rank $|H|$.
%
%
%
%
Thus, there exists a set $U:=\{u^1, \dotsc, u^t\} \subseteq \R^{|H|} \setminus\{0\}$ that generates the extreme rays of $K$, i.e.,
%
\[
K = \bigg\{\sum_{i=1}^t \lambda_i u^i : \lambda_i \ge 0 ~ \forall ~ i \in \{1, \dotsc, t\}\bigg\}
\]
and $u^i$ satisfies $|H|-1$ linearly independent constraints in~\eqref{eqCone} at equality for each $i \in \{1, \dotsc, t\}$.
%

\begin{CLAIM}\label{claimExtRays}
For each $\tilde{u} \in U$ we have
\begin{equation}\label{eqExtRaySupp}
|\supp(\tilde{u})| \le m+1 .
\end{equation}
Also, each $\tilde{u}\in U$ can be scaled to have integer components and satisfy $\|\tilde{u}\|_{\infty}\le \Delta$.
\end{CLAIM}
\begin{proof}
Set $T:=\supp(\tilde{u})$ and without loss of generality assume $T =\{ 1,\dots ,|T|\}$.
Recall that $\tilde{u}$ satisfies a set of $|H|-1$ linearly independent constraints in~\eqref{eqCone} at equality.
One such set is composed of $|H|-|T|$ constraints from the system $D_1 \tilde{u} \le 0, ~ D_2 \tilde{u} \ge 0$ and $|T|-1$ constraints from $0=A_H\tilde{u}=A_T\tilde{u}_T$.
By this choice of constraints, it follows that
\[
|\supp(\tilde{u})| = |T| \le m+1
\]
and $|T| - 1 \le \rank(A_T) \le \min\{|T| ,m \}$.
Recalling $n > m$ and $\rank(A) = m$, the latter inequalities imply that there exists an index set $\overline{T}$ satisfying $T \subseteq \overline{T} \subseteq \{1, \dotsc, n\}$, $A_{\overline{T}} \in  \Z^{m\times (m+1)}$ and $\rank (A_{\overline{T}}) = m$.

Let $\bar{u} \in \R^{m+1}$ denote the vector obtained by appending $m+1-|T|$ zeros to $\tilde{u}_T$.
There exists an index set $I \subseteq \overline{T}$ with $|I| = m$ and $A_I$ invertible.
Let $i$ denote the singleton in $\overline{T} \setminus I$.
Because $A_{\overline{T}} \bar{u} = 0$, we have
\[
A_I \bar{u}_I = - A_i \bar{u}_i.
\]
If $\bar{u}_i = 0$, then $\bar{u} = 0$ and so $\tilde{u} = 0$.
However, this contradicts that $\tilde{u} \in U$.
Hence, $\bar{u}_i \neq 0$.
Scale $\bar{u}$ such that $|\bar{u}_i|=|\det ({A}_I)|$.
Applying Cramer's rule demonstrates that
\[
|\bar{u}_j| = |\det(A_{I \cup \{i\} \setminus \{j\}})| \qquad \forall ~ j \in I.
\]
Hence, $\bar{u}$, and consequently $\tilde{u}$, can be scaled to have integer components with $\|\tilde{u} \|_\infty  \leq\Delta$.
This proves the claim.
\qed
\end{proof}
For the rest of the proof, we assume that each $\tilde{u} \in U$ is scaled such that the conclusions of Claim~\ref{claimExtRays} hold.

Recall $\bar{z}_H -x^*_H \in K$.
By Carath\'eodory's theorem, there exists an index set $I\subseteq \{ 1,\dots , t \}$ with $|I| =|H|$ and $\lambda_i \in \mathbb{R}_{\ge 0}$ for each $i\in I$ such that

%
\[
  \bar{z}_H -x^*_H=\sum_{i\in I } \lambda_i u^{i} .
\]
Set $w := {\sum_{i\in I}} \lfloor\lambda_i \rfloor u^{i}$.
Using standard techniques in proximity proofs, it can be verified that $\tilde{z}:= \bar{z}_H - w
$ is a feasible solution to \eqref{IP2}, and $\tilde{x}:= x^*_H+ w
$ is a feasible solution to \eqref{LP2} (see, e.g.,~\cite[Theorem 1]{CooSens}).
Moreover, 
\begin{equation}\label{eqProjectedObjective}
 c^\transpose_H \bar{z}_H+c^\transpose_H x^*_H= c^\transpose_H \tilde{z}+ c^\transpose_H \tilde{x} .
\end{equation}
Because $x^*_H$ is optimal for~\eqref{LP2}, we have $c_H^\transpose\tilde{x}\le c_H^\transpose x^*_H$.
Combining this with~\eqref{eqProjectedObjective} proves that $c_H^\transpose \tilde{z}\ge c_H^\transpose \bar{z}_H$.
Because $\bar{z}_H$ is optimal for~\eqref{IP2}, we have $c_H^\transpose \tilde{z}= c_H^\transpose \bar{z}_H$ and  $\tilde{z}$ is also an optimal solution to \eqref{IP2}.

Define $z^* \in \Z^n_{\ge 0}$ component-wise to be
\[
z^*_i := \begin{cases} \tilde{z}_i &\text{if } i \in \{1, \dotsc,  |H| \}\\ 0 &\text{otherwise}.\end{cases}
\]
By construction $z^*$ is an optimal solution to~\eqref{IP} because $Az^* = A_H \tilde{z} = A_H \bar{z}_H = b$ and $c^\transpose z^* = c_H^\transpose \tilde{z}=c_H^\transpose \bar{z}_H=c^\transpose \bar{z}$.
%
%
By~\eqref{eqDefnH} and~\eqref{eqExtRaySupp}, we arrive at the final result:
\begin{align*}
 \| z^*  -x^* \|_1 =\| \tilde{z} -x^*_H \|_1
&\le {\sum_{i\in I}} (\lambda_i - \lfloor \lambda_i \rfloor) \|u^{i}\|_1 \\
&< {|H|} \cdot  (m+1)\cdot \|u^i\|_{\infty} \\
&\leq (m+1)\cdot (m+S)\cdot \Delta.\qedEQN
\end{align*}
%
%
\end{proof}

\begin{proof}[of Theorem~\ref{DeltaSpars}]
Let $z^*$ be an optimal solution of~\eqref{IP} with minimum support.
The definition of $S$ states that ${S} := |\supp (z^* )|$.
Define $\tilde{A}\in\mathbb{Z}^{m\times {S}}$ as the submatrix of $A$ corresponding to the support of $z^*$.
If ${S}\leq 2m$, then the result holds.
Thus, assume that ${S}> 2m$ which implies $\log( \sfrac{{S}}{m})<(\sfrac{{S}}{m})-1$.
Theorem 1 in Aliev et al.~\cite{ADEOW2018} states that
\[
{S} < m+\log \left(\sqrt{\det (\tilde{A}{\tilde{A}}\strut^{\transpose})}\right).
\]
The Cauchy-Binet formula for $\det (\tilde{A}\tilde{A}\strut^\transpose)$ states that
\[
\det(\tilde{A} \tilde{A}\strut^\transpose) = \sum_{\substack{B \text{ is an } m \times m \\\text{submatrix of } \tilde{A}}} \det(B)^2.
\]
See, e.g.,~\cite{HJ2012}.
Combining the previous inequalities yields
\begin{align*}
{S} < m+\log \left(\sqrt{\det (\tilde{A}\tilde{A}\strut^\transpose)}\right)&  \leq m+\log \left(\sqrt{\binom{{S}}{m}\Delta^2}\right)\\[.1 cm]
&\leq m+\log \left({S}^{m/2}\Delta \right)\\
&=m+\frac{m}{2}\log\left( \frac{{S}}{m}\right)+\frac{m}{2}\log(m) +\log (\Delta) \\&
< m+\frac{m}{2}\left(\frac{{S}}{m}-1\right)+\frac{m}{2}\log (m) +\log (\Delta) \\&
=\frac{{S}}{2}+\frac{m}{2}+\frac{m}{2}\log (m) +\log (\Delta).
\end{align*}
Therefore,
\[
|\supp(z^*)| = {S}<m+m\log (m) +2\log (\Delta) \leq 2m\log \bigl(\sqrt{2m} \cdot \Delta^{1/m}\bigr).\qedEQN
\]
%

%
\end{proof}

\begin{proof}[of Corollary \ref{DeltaSparsMIP}]
  Let $z^*$ be {an optimal \eqref{MIP} solution} with minimal support.
  By applying Theorem \ref{DeltaSpars} to the standard form integer program with constraint matrix $A_{\mathcal{I}}$ and right hand side $b-A_{\mathcal{J}}z^*_{\mathcal{J}}{\in\mathbb{Z}^m}$, where $\mathcal{J}:=\{1,\dots,n\}\setminus \mathcal{I}$, we see that
  \[
  |\supp (z^*_{\mathcal{I}} )|< 2m\log \bigl(\sqrt{2m}\cdot \Delta^{1/m}\bigr).
  \]
    Similarly, by considering the standard form linear program with constraint matrix $A_{\mathcal{J}}$ and right hand side $b-A_{\mathcal{I}}z^*_{\mathcal{I}}$, we see that
  \(
|\supp(z^*_{\mathcal{J}})|\le m.
  \)
  Hence,
  \[
  |\supp(z^*)| \le  |\supp(z^*_{\mathcal{J}})| + |\supp(z^*_{\mathcal{I}})|   \le m + 2m\log \bigl(\sqrt{2m}\cdot \Delta^{1/m}\bigr).\qedEQN
  \]
\end{proof}

\section{Results on proximity}\label{sec:proximity}

%
\begin{proof}[of Theorem \ref{DeltaProx}]
For now, assume that $m \ge 2$.
Combining Lemma \ref{MixedProx} with Theorem~\ref{DeltaSpars} demonstrates that there exists an optimal \eqref{IP} solution $z^*$ satisfying
\[
\begin{split}
\|z^* -x^*\|_1 &<(m+1)\cdot (m+S)\cdot \Delta \\
&\leq (m+1) \cdot 2m\log (2\sqrt{m}\cdot \Delta^{1/m})\cdot\Delta\\
&\leq 3m^2 \log (2\sqrt{m}\cdot \Delta^{1/m})\cdot\Delta.
\end{split}
\]
This completes the proof when $m \ge 2$.

It is left to consider the case $m  = 1$.
Here we have $\Delta = \|A\|_{\infty}$.
The proximity bound~\eqref{eqEWEntryBound} states that there exists an optimal solution $z^*$ to~\eqref{IP} satisfying
\[
\|z^* - x^*\|_1 \le 2 \cdot \|A\|_{\infty}+1 = 2\Delta+1 < 3\Delta \le 3m^2 \log (2\sqrt{m}\cdot\Delta^{1/m})\cdot \Delta .
\]
This completes the proof.
\qed
\end{proof}


We use the following result to prove Theorem~\ref{eqThmEntry}.
\begin{lem}[Theorem 1 in~\cite{SumSubdet}]\label{lemAppx}
For every $\epsilon > 0$, there exists an $m \times m$ submatrix $B$ of $A$ satisfying
\[
\Delta \le |\det(B)| \cdot (e \cdot \ln\left((1+\epsilon) \cdot m)\right)^{m/2}.
\]
The matrix $B$ can be found in time that is polynomial in $m,n$, and $1/\epsilon$.
\end{lem}
Notice that selecting $\epsilon = 1/m$ in Lemma~\ref{lemAppx} yields the approximation factor
\begin{align*}
\Delta & \le |\det(B)| \cdot \left(e \cdot \ln\left(m+1\right)\right)^{m/2} \\
& = |\det(B)| \cdot  \left(\frac{e}{\log(e)} \cdot \log(m+1)\right)^{m/2}\\
& \le |\det(B)| \cdot (2 \cdot \log(m+1))^{m/2},
\end{align*}
which is precisely~\eqref{eqDetApprox}.

\begin{proof}[of Theorem~\ref{eqThmEntry}]
Let $B$ be an $m\times m$ submatrix of $A$ satisfying~\eqref{eqDetApprox}.
There exists a unimodular matrix $U \in \Z^{m\times m}$ such that $UB$ is an upper triangular matrix with non-negative diagonal entries $d_i$.
The so-called row-style Hermite normal form is a transformation satisfying these properties, and the corresponding unimodular matrix $U$ can be computed in polynomial time, see~\cite{AS1986}.

Unimodular matrices preserve the absolute value of $m\times m$ determinants, so $|\det(B)| = |\det(UB)|$.
By~\eqref{eqDetApprox} we see that
\[
\begin{array}{rl}
\displaystyle \Delta \le |\det(B)| \cdot  (2 \cdot \log(m+1))^{m/2} & =\displaystyle   |\det(UB)| \cdot  (2 \cdot \log(m+1))^{m/2}\\[.25 cm]
&\displaystyle = \bigg(\prod_{i=1}^m d_i \bigg)\cdot  (2 \cdot \log(m+1))^{m/2} \\ [.5 cm]
&\displaystyle \le \|UB\|_{\infty}^m \cdot  (2 \cdot \log(m+1))^{m/2}.
\end{array}
\]

Therefore, we may apply Theorem~\ref{DeltaProx} to obtain the bound
\begin{align*}
\|z^* -x^*\|_1  < &3m^2 \log (2\sqrt{m}\cdot \Delta^{1/m}) \cdot \Delta\\
 < & 3m^2  \log (2\sqrt{2m\log (m+1)} \|UB\|_\infty  ) \cdot (2\log (m+1))^{m/2} \cdot\|UB\|_\infty ^m.
\end{align*}
This completes the proof.
\qed
\end{proof}

Next, we present a proof of Corollary \ref{corProxPlusSparsity}.
We advise the reader that the proof is similar to the proof of Lemma~\ref{MixedProx}.

\begin{proof}[of \rm{Corollary \ref{corProxPlusSparsity}}]
Let $z^*$ be an optimal solution to \eqref{GIP} with minimal support on the first $n$ components.
Consider the $n$-dimensional integer program obtained by fixing the last $d$ variables of \eqref{GIP} to the last $d$ components of $z^*$.
Similarly, consider the $n$-dimensional linear program obtained by fixing the last $d$ variables of (GLP) to the last $d$ components of $x^*$.
These lower dimensional problems are in standard form.
Therefore, by applying Theorem \ref{DeltaSpars} to these lower dimensional problems and recalling $\Delta_m(A) \le \delta$, we see that
\[
\bar{S}:=|\supp (x^* ) \cup \supp (z^* )|\leq \min \bigl\{ n,m+2m\log (\sqrt{2m}\cdot\delta^{1/m}) \bigr\} +d .
\]
The conclusion of Corollary \ref{corProxPlusSparsity} follows directly from the inequality
\[
\|z^*-x^*\|_1 <\min\{m+t+1,n+d\} \cdot \bar{S} \cdot \delta,
\]
and we complete the proof of the corollary by demonstrating that the latter inequality holds.
For this, we project the original optimization problems onto the union of the supports of $x^*$ and $z^*$.
To simplify the presentation, we assume that $\bar{S}=n+d$.
The other case is similar.

As in the proof of Lemma~\ref{MixedProx}, we create a pointed cone from the constraints defining~\eqref{GIP}.
In order to assure a pointed cone, we introduce redundant constraints.
Let $b_3 \in \Z^d$ be the vector where every component is $\lceil \|x^*\|_\infty \rceil+ \|z^*\|_\infty$, and define
\[
D :=
\begin{pmatrix}
 0 & C \\-I_n & 0 \\ 0 & I_d
 \end{pmatrix}\in\mathbb{Z}^{(t+n+d)\times (n+d )} \quad \text{and} \quad
 f=\begin{pmatrix}
 b_2 \\ 0 \\ b_3
 \end{pmatrix}\in\mathbb{Z}^{t+n+d}.
\]
By construction, $z^*$ is an optimal solution to the integer program
 \[
 \max\left\{ c^\transpose z  : [A,B] z = b_1, ~Dz \leq f,~z \in\mathbb{Z}^{n+d} \right\},
 \]
and $x^*$ is an optimal vertex solution to the corresponding linear relaxation.

 Subdivide the rows of $D$ such that $D_1 z^* < D_1 x^* $ and $D_2 z^* \geq D_2 x^*$.
 Define the polyhedral cone
 %
%
 %
 \[
 K :=\left\{ u\in\mathbb{R}^{n+d } : [A,B]u=0,~\begin{pmatrix} D_1 \\ - {D}_2 \end{pmatrix} u \le 0\right\}.
 \]
 Observe that $z^*-x^*\in K$.
 Moreover, the introduction of $b_3$ and $I_d$ ensures that $\rank (D)=n+d$ and that $K$ is pointed.

 Let $U$ be a finite set generating the extreme rays of $K$.
The proof of Claim~\ref{claimExtRays} can be used to demonstrate that for each $\tilde{u} \in U$ we have
 \[
 |\supp(\tilde{u})|\le \min\{m+t+1,n+d\},
 \]
 and each $\tilde{u} \in U$ can be scaled such that $\tilde{u}\in\mathbb{Z}^{n+d}$ and $\|\tilde{u}\|_{\infty}\le \delta$.
The main difference that arises when repeating this proof is that the matrices $D_1$ and $D_2$ contain rows of $[0, C]$ rather than simply rows from the identity matrix as was the case in the proof of Lemma~\ref{MixedProx}.
It is this difference that necessitates the choice of the data parameter $\delta$ and dictates its definition.

By Carath\'eodory's theorem, there exist $k \le \bar{S}$ vectors $u^{1},\dots ,u^{k} \in U$ and coefficients $\lambda_1 ,\dots ,\lambda_k \in \mathbb{R}_{\ge 0}$ such that
 $
 z^* -x^*=\sum_{i=1}^k \lambda_i u^{i} .
 $
 Because $z^* -\sum_{i=1}^k \lfloor\lambda_i \rfloor u^{i}$ is also an optimal solution to \eqref{GIP}, we can assume, without loss of generality, that $\lambda_1, \dotsc, \lambda_k < 1$ (the reasoning is similar to that in the proof of Lemma~\ref{MixedProx}).
 This implies that
 \begin{equation*}
 \| z^* -x^* \|_1 \le \sum_{i=1}^k \lambda_i \|u^{i}\|_1 <\min \{ m+t+1,n+d \}\cdot \bar{S}\cdot \delta.\qedEQN
 \end{equation*}
 \end{proof}

\noindent\textbf{Acknowledgements.} J. Lee was supported in part by ONR grant N00014-17-1-2296.

\bibliographystyle{splncs04}
\bibliography{references}
\end{document}